\documentclass[11pt]{amsart}

\usepackage{amssymb} 
\usepackage[latin1]{inputenc}
\usepackage[matrix,arrow,curve]{xy}
\usepackage[mathscr]{eucal}
\usepackage{epic,eepic}
\usepackage{bbm}

\makeindex
\newtheorem{theorem}{Theorem}[section]
\newtheorem{exptheorem}[theorem]{Expected Theorem}
\newtheorem{lemma}[theorem]{Lemma}
\newtheorem{prop}[theorem]{Proposition}

\theoremstyle{definition}

\newenvironment{pf}
{\medskip\noindent {\it Proof --- \ }}
{\hfill\nobreak $\Box$ \par\bigbreak}

\newcommand{\RR}{{\mathcal R}}

\newcommand{\Hom}{\text{Hom}}

\newcommand{\Qp}{{\mathbb Q}_p }
 
\newcommand{\C}{{ \mathbb C  }}
\newcommand{\R}{{ \mathbb R  }}
\newcommand{\Q}{{ \mathbb Q } }
\newcommand{\Z}{{ \mathbb Z  }}

\newcommand{\Lie}{\text{Lie}}
\newcommand{\m}{m}

\newcommand{\Gal}{{\mathrm{Gal}\,}}

\newcommand{\A}{\mathbb A}

\newcommand{\anneau}{{ \mathcal O}}

\newcommand{\End}{{\text{End}}}

\newcommand{\U}{{\text{U}}}

\newcommand{\Gl}{{\text {GL}}}

\newcommand{\GL}{{\text {GL}}}

\newcommand{\spec}{{\text{Spec\,}}}

\newcommand{\G}{{\mathfrak g}}

\newcommand{\ses}{{\text{ss}}}

\newcommand{\Frob}{{\text{Frob\,}}}

\newcommand{\rhob}{{\bar \rho}}

\renewcommand{\v}{{\bf v}}

\newcommand{\Qpb}{{\overline{\Q}_p}}
\renewcommand{\Gal}{{\rm Gal}}

\newcommand{\OO}{\mathcal{O}}

\newcommand{\Qb}{\overline{\Q}}
\newcommand{\AAA}{\mathbb{A}}

\renewcommand{\G}{G}

\newcommand{\crys}{\text{crys}}

\newcommand{\Gt}{\tilde{G}}

\begin{document}

\baselineskip 14.8pt

\bibliographystyle{style} 
\title[The sign of Galois representations]{The sign of Galois representations attached to automorphic forms for unitary groups}
\date{April 2, 2008}
\thanks{We thank Laurent Clozel, Michael Harris and Jean-Pierre Labesse, for many useful conversations and for their constant support. This paper
rely on the book project \cite{book}, and we thank all its authors for having made it possible. During the elaboration and writing 
of this paper, Jo\"el Bella\"iche was supported by the NSF grant DMS 05-01023, and Ga\"etan Chenevier was supported by the C.N.R.S}
\author[J.~Bella\"iche]{Jo\"el Bella\"iche}
\email{jbellaic@brandeis.edu}
\address{Jo\"el Bella\"iche\\Brandeis University\\
415 South Street\\Waltham, MA 02454-9110\\U.S.A}
\author[G. Chenevier]{Ga\"etan Chenevier}
\email{chenevie@math.univ-paris13.fr}
\address{Ga\"etan Chenevier\\C.N.R.S, L.A.G.A., Universit\'e Paris 13\\ 
99 Av. J-B. Cl\'ement\\93430 Villetaneuse\\France}

\maketitle


\section{Introduction}

\subsection{The sign of a representation}

Let $L$ be a field of characteristic $0$ or  greater than $2$.
Let $G$ be a group, and $g \mapsto g^c$ an involution of $G$. 
For $\rho$ a representation $G \rightarrow \Gl_n(L)$,
we define $\rho^\bot : G \rightarrow \Gl_n(L),\ g \mapsto {}^t\rho(g^c)^{-1}$.
The equivalence class of the representation $\rho^\bot$ only depends on the 
equivalence class of $\rho$.

We fix $\chi:G \rightarrow L^\ast$ a character such that $\chi(g)=\chi(g^c)$
for all $g$. This property ensures that 
$\rho \mapsto \rho^\bot \chi^{-1}$ is an involution. 
In the applications, $G$ will be the absolute Galois group of a CM field $K$,
$c$ the outer automorphism defined by the non trivial element in $\Gal(K/F)$
where $F$ is the maximal totally real field in $K$, and $\chi$ will be
a power of the cyclotomic character.  

Let $\rho$ be a semi-simple representation $G \rightarrow \Gl_n(L)$ such that
\begin{eqnarray}\label{autobot} \rho^\bot \simeq \rho \chi. \end{eqnarray}
This property is obviously stable by extension of the field of coefficients $L$.

We shall now attach to any {\it absolutely irreducible} 
$\rho$ satisfying (\ref{autobot}) an invariant, that we call its {\it sign}. The invariant can take the value $+1$ or $-1$.
By Schur's lemma there exists a unique (up to a scalar) matrix 
$A \in \Gl_n(L)$ such that 
\begin{eqnarray} \label{A} \rho^\bot = A \rho A^{-1} \chi.\end{eqnarray}
 Applying this relation twice, we see that $A {}^t A^{-1}$ commutes with $\rho^\bot$, hence by Schur's lemma again is a scalar matrix $\lambda$. So ${}^t A = \lambda A$ and $\lambda = \pm 1$. This sign is called the {\it sign} of {$\rho$}. Note that it is necessarily $1$ if $n$ is odd, since there is no invertible antisymmetric matrix in odd dimension.

It is obvious that the sign of an absolutely irreducible representation
is unchanged by arbitrary extensions of the coefficient field $L$.

\subsection{The context: the book project on Galois representations attached to unitary groups} 

Let $F$ be a totally real field and $K$ a totally imaginary quadratic 
extension, $c \in \Gal(K/F)$ the non-trivial automorphism.
Let $\Pi$ be a cuspidal automorphic representation for $\Gl_n$ over $K$, and assume that $\Pi$ is polarized (i.e. the contragredient $\Pi^\vee$ of $\Pi$ is isomorphic to $\Pi \circ c$) and that $\Pi_\infty$ is algebraic regular (see \cite[General Hypotheses 2.1]{harris}).  

Under those hypotheses, the many  coauthors of the {\it book project} (\cite{book})
expect to prove the existence of an attached compatible system of 
Galois representations (see \cite[Theorem 4.2.4]{harris}):
\begin{exptheorem}[Book Project] \label{galois} There is a number field $E(\Pi)$ 
and a compatible system $\rho_{\Pi,\lambda} : 
G_K \rightarrow \GL(n,E(\Pi)_\lambda)$ of $\lambda$-adic representations, where $\lambda$ runs through finite places of $E(\Pi)$ such that
for all finite primes $v$ of $K$ of residue characteristic prime to $N_{E(\Pi)/\Q}(\lambda)$,
$$\rho_{\Pi,\lambda}^{F-ss}|_{\G_v} \simeq {\rm{L}}(\Pi_v \otimes | \bullet |_v^{\frac{1-n}{2}}),$$
where $G_v$ is a decomposition group of $K$ at $v$ and $L(\bullet)$ is the local Langlands
correspondence. 
\end{exptheorem}
The given property suffices to characterize uniquely
$\rho_{\Pi,\lambda}$ up to isomorphism and implies that $\rho_{\pi,\lambda}$ satisfies (\ref{autobot}); more precisely, let $c$ be a complex conjugation in $K$, that is  an element of $G_F-G_K$ of order $2$. We set $g^c = cgc^{-1} = cgc$
for $g \in G_K$: this is an automorphism of order $2$. For that automorphism,
we have
$$\rho^\bot(g) = {}^t \rho (g^c)^{-1} \simeq \rho(g) \omega(g)^{n-1}$$
where $\omega$ is the cyclotomic character.

The proof of this theorem relies on a special case, which is the same theorem with stronger hypotheses on the extension $K/F$ and on the automorphic representation $\Pi$ : see \cite[Expected Theorem 1.4]{harris}. This theorem is the output of the comparison of the two stabilized trace formulas, and other works,
to be done in books one and two of the book project. Book one has been 
mainly written and most chapters are available on \cite{book}, while the writers of book two are expected to handle their chapters by May 2008. 

The derivation from \cite[expected Theorem 1.4]{harris} to the theorem quoted above is carefully written in \cite{harris}. We shall use several lemmas proved by Harris during this this derivation.

The theorem also includes many other specifications on $\rho_{\Pi,\lambda}$, including the expected determination on the Hodge-Tate weights of $\rho_{\Pi,\lambda}$ at places of same residual characteristic as $\lambda$ (see also \S\ref{generalnotations} below). This description implies, since $\Pi_\infty$ is cohomological, that these weights are distinct integers, hence that $\rho_{\Pi,\lambda}$ is a direct sum of non-isomorphic absolutely irreducible representations of $G_K$.

\subsection{The result} The object of this article is to prove 
(admitting \cite[Expected Theorem 1.4]{harris})
\begin{theorem}\label{main}
For every finite prime $\lambda$ of $E(\Pi)$, every irreducible
 factor of $\rho_{\Pi,\lambda}$ that 
 satisfies (\ref{autobot}) has sign $+1$.
\end{theorem}

It is expected that $\rho_{\Pi,\lambda}$ is absolutely irreducible (this is known if $n \leq 3$ by \cite{BR} and in many cases if $n=4$ by an unpublished work of Ramakrishnan). 
If it is so, $\rho_{\Pi,\lambda}$ has only one factor
and satisfies (\ref{autobot}), and our theorem simply asserts that 
its sign is $+1$: this is obvious when $n$ is odd, but new when $n$ is even.

\subsection{Historical remarks}
 
The question of the sign of Galois representations attached to 
polarized automorphic representations of $\Gl_n$ on a totally real or CM field
is out at least since Clozel, building on the work of Kottwitz, proved
their existence in many cases in the mid nineties. 
More recently, this question has been extensively discussed in
\cite{CHT} where some cases of the above theorem, concerning 
Galois representations with some constraining properties ensuring
they have a nice and workable deformation theory, are proved by a 
very indirect method -- indeed the whole long and hard paper is written
with an unknown sign $\epsilon$ and only near the end, after the Taylor-Wiles method has been adapted to unitary groups, is it shown that $\epsilon=-1$ 
leads to a contradiction! 

\par \bigskip

Theorem~\ref{main} appears, without its proof, in the concluding remarks of our book \cite{BC} (see \cite[theorem 9.5.1]{BC}) that was made public late 2006.
We knew the proof that follows then\footnote{At least for places $\lambda$ of residual characteristic $p$ split in $K$.}, and told it to a few colleagues, but 
decided to wait for a more advanced version of the book project, on which it 
depends, before writing it.

\par \bigskip
Meanwhile, one of us, Ga\"etan Chenevier, together with Laurent Clozel,
has found a completely different proof of a special case of 
Theorem~\ref{main}: they prove that the Galois representation attached to a polarized and cohomological automorphic representation of $\Gl_n$ over a totally real field, which is square integrable at some 
finite place, is symplectic\footnote{It may seem strange at first glance that
this symplecticity result is a special case of our sign $+1$ result -- one would naively expect $-1$. Actually, a closer look shows that the two results are consistent. Let us explain why. 

Chenevier and Clozel's result is about an even dimensional absolutely irreducible representation of 
$r$ of $G_F$, say of dimension $2n$, $r: G_F \rightarrow \GL_{2n}(L)$, as constructed by Clozel (the irreducibility following from works of Harris-Taylor and Taylor-Yoshida). The representation $r$ satisfies $r^\vee \simeq r\omega^{2n-1}$ so
$${}^t r(g)^{-1} = P r(g) P^{-1} \omega(g)^{2n-1},\, \,\, \,\forall\, g\, \in\, G_F,$$
where $\omega$ is again the cyclotomic character, and Chenevier-Clozel result is that $P$ is anti-symmetric. We claim that this implies that the restriction $\rho = r|G_K$, that satisfies clearly 
(\ref{autobot}), has sign $+1$.

Indeed, by changing the basis, we may assume that $r(c)$ is diagonal, hence symmetric. Moreover, the relation of selfduality of $r$ above implies
$$ r(c) = - P r(c) P^{-1}$$
since $\omega(c)=-1$. An immediate computation shows then that ${}^t \rho(g^c)^{-1} = A \rho(g) A^{-1} \omega(g)^{2n-1}$ with $A = r(c) P$. But we see that
$${}^t A = {}^t P {}^t r(c) = (-P) r(c) = r(c) P = A$$
hence the sign of $\rho$ is $+1$.}
when $n$ is even. Their proof (\cite{CC}) is less 
expensive in difficult tools than ours, using simply the new insight in the trace formula they discovered in \cite{CC}. It does not seem that their method can be extended to the case of a CM field, or even to the case of an automorphic representation that does not satisfy any local square-integrability hypothesis. 

\par \bigskip
Let us mention also that in a recent preprint \cite{gross}, B. Gross introduces a general notion 
of {\it odd} Galois representations and conjectures that the expected Galois representations 
attached to definite reductive groups $G$ are odd in his sense. Our theorem proves his conjecture 
when $G$ is the a unitary group attached to a CM extension $K/F$, in which case it has the following 
meaning.

Let $\Gt$ be the semi-direct product of $\Gal(K/F)=\langle c \rangle = \Z/2\Z$ by $\GL_n(L) \times L^*$ 
with respect to the order two automorphism $(x,y) \mapsto (y {}^tx^{-1},y)$ (see \cite[Ch.I]{CHT} for similar considerations). 
Assume that $\rho : G_K \rightarrow \Gl_n(L)$ satisfies ($\ref{autobot}$), is absolutely irreducible, and fix $A$ a matrix as in (\ref{A}),
and $\epsilon = \pm 1$ the sign of $\rho$. Consider the morphism $G_K \rightarrow \GL_n(L) \times L^*$ 
defined by $g \mapsto (\rho(g),\chi(g)^{-1})$. A simple computation shows that this map extends to a 
morphism $\tilde{\rho} : G_F \rightarrow \Gt$ if we set $\tilde{\rho}(c)=( {}^t\!A^{-1},\epsilon) c$. 
Assume now that $\rho=\rho_{\Pi,\lambda}$. The map $\tilde \rho$ is the analog in our situation of the 
map denoted $\rho$ whose existence is conjectured in \cite[page 8]{gross} and 
Gross predicts that the conjugation by $\tilde{\rho}(c)$ on $\Lie(\Gl_n)$ is a Cartan involution, that is, has the form  
$X \mapsto - P {}^t \!X P^{-1}$ with $P$ a {\it symmetric} invertible matrix. 
In our situation, the conjugation by $\tilde{\rho}(c)$ on the Lie algebra is 
the map $X \mapsto - {}^t\!A^{-1} {}^t X {}^t\!A$. So we see that Gross' 
prediction amounts to ``$A$ is symmetric'', which is exactly our theorem.
\footnote{When $\rho$ is not assumed to be irreducible anymore, note that 
theorem~\ref{main} still implies that we may find some symmetric $A$ such 
that (\ref{A}) holds, hence a $\tilde{\rho}$ as above satisfying Gross' 
conjecture.}

\subsection{Idea of the proof}

The idea of the proof is very simple. Assume that we know that the representation
$\rho_{\Pi,\lambda}$ is irreducible. Then there is nothing to prove if $n$ is odd.
When $n$ is even, we can reduce to the odd case, as follows: descend $\Pi$ to a unitary group in $n$ variables, transfer the result to an automorphic representation $\pi$ 
of a unitary groups in $n+1$ variables which is compact at infinity, using a special case of endoscopic transfer proved by Clozel, Harris and Labesse. Use eigenvarieties to deform $\pi$ 
into a family of automorphic forms whose Galois representations are generically irreducible. For  those Galois representations, the sign is $+1$ since their dimension is odd. Specialize this result to deduce that the components of the
representation attached to $\pi$, including $\rho_{\Pi,\lambda}$, have sign $+1$.

There are several technical difficulties that make the proof a little bit more indirect:
in the current state of science, we do not know that $\rho_{\Pi,\lambda}$ is 
(absolutely) irreducible,
and we cannot descend $\Pi$ to $\U(n)$ or transfer it to $\U(n+1)$ without supplementary assumptions on $K/F$ and $\Pi$. Moreover,  we cannot always deform a representation $\pi$ in a family whose Galois representation is generically irreducible. But this is not a big issue, since, as was already observed 
in \cite[\S7.7]{BC},
we can actually do so in two steps, deforming $\pi$ in a family whose generic members can themselves be deformed irreducibly. Similarly the obstacle posed by the conditions on descent and endoscopic transfer can be solved by base change techniques inspired by the ones used in \cite{harris}.

\subsection{Notations and conventions}\label{generalnotations} Our general convention will be that the local Langlands correspondence is normalized so that geometric Frobeniuses correspond to uniformizers (and as in \cite{HT}). If $\pi$ is an unramified complex representation of $\GL_n(E)$ with $E$ a $p$-adic local field, or more generally an irreducible smooth representation with a nontrivial vector fixed by a Iwahori subgroup, we shall often denote by $L(\pi)$ the semisimple conjugacy class in $\GL_n(\C)$ of the geometric Frobenius in the $L$-parameter of $\pi$.

If $K$ is a field, we shall denote by $G_K$ its absolute Galois group 
$\Gal(\overline{K}/K)$; when $K$ is a number field and $v$ a place of $K$ we also write $G_v$ for $G_{K_v}$.
	
	We shall use the following notions of $p$-adic Hodge theory. Let us fix $E$ a finite extension of $\Q_p$, $\Qpb$ an algebraic closure of $\Qp$, and let $V$ be a $p$-adic representation of $G_E$ of dimension $n$ over $\Qpb$. To such a representation Sen attaches a monic polynomial $P_{\rm sen}(T) \in (\Qpb\otimes_{\Q_p} E)[T]$ of degree $n$, whose roots will be called the Hodge-Tate weights of $V$ (even when they are not natural integers). Our normalization of the Sen polynomial is the one such that the Hodge-Tate weight of the cyclotomic character $\Qpb(1)$ is $-1 \in \Qpb\otimes_{\Q_p} E$. Under the natural identification $\Qpb\otimes_{\Q_p} E = \Qpb^{\Hom(E,\Qpb)}$, we shall often write them as a collection $\{k_{i,\sigma}\}$ for all $i \in \{1,...,n\}$ and all $\sigma \in \Hom(E,\Qpb)$, ordered so that for each embedding $\sigma$ we have $$k_{1,\sigma} \leq k_{2,\sigma} \leq \dots \leq k_{n,\sigma}.$$ 
	
	We shall need to consider various partial sums of those
weights, for which the following definitions will be useful. For
$I$ a subset of $\{1,\dots,n\} \times \Hom(K_w,\Qpb)$, we denote by 
$k_I$ the sum  $\sum_{(i,\sigma) \in I} k_{i,\sigma}.$
When $I=\{i\} \times \Hom(K_w,\Qpb)$, we write $k_i$ instead of $k_I$.
Thus $k_i = \sum_\sigma k_{i,\sigma}$. 

		Assume now that $V$ is crystalline in the sense of Fontaine. Let $E_0 \subset E$ be the maximal unramified extension of $\Q_p$ inside $E$, and let $\v : \Qpb \rightarrow \Q$ be the valuation {\it normalized so that $\v(p)=e$}, where $e$ is the absolute ramification index of $E$.	Fontaine attaches to $V$ an $E_0$-vector space $D_{\rm crys}(V)$ with a semilinear action of the crystalline Frobenius $\varphi$ (commuting with $\Qpb$), and which is free of rank $n$ over $E_0\otimes_{\Q_p} \Qpb$. If $f=[E_0 : \Q_p]$, then $\varphi^f$ is $\Qpb\otimes_{\Q_p} E_0$-linear and commutes with $\varphi$, so its characteristic polynomial $P_\varphi(T)$ actually belongs to $\Qpb[T]$.
 This polynomial will be referred as the {\it characteristic polynomial} of $\varphi$, its roots are the {\it eigenvalues} of the crystalline Frobenius, and their valuations (with respect to $\v$) its {\it slopes}.\footnote{This definition is slightly different from the usual definition of the slopes of an isocrystal (which are ours divided by $[E:\Q_p]$), but it will be convenient to us.} With these notations, if the $k_{i,\sigma}$ are the Hodge-Tate weights of $V$, then the weak admissibility property of $D_{\rm crys}(V)$ implies in particular that $$\v(P_\varphi(0)) = \sum_{i,\sigma} k_{i,\sigma}.$$

	We can now explain a bit more precisely the $p$-adic part the Expected theorem~\ref{galois}. Assume that $w$ is a finite place of $K$ with the same residual characteristic as $\lambda$, and assume than $\Pi_w$ is unramified. Let $P_w(T) \in E(\Pi)[T]$ be the characteristic polynomial of $L(\Pi_w|.|^{(1-n)/2})$. Then part of Expected Theorem~\ref{galois} asserts that {\it $\rho_{\Pi,\lambda}|G_w$ is a crystalline representation} (see \cite{harris}). 
	
	At present, it is not clear whether we will be able to identify the characteristic polynomial $P_\varphi \in E(\Pi)_\lambda[T]$ of the crystalline Frobenius with the image of $P_w(T)$ in this full generality, as it should be. However, we will know that $P_\varphi=P_w$ by construction under the following extra assumptions : assumptions (H1) and (H2) stated in \S\ref{proofone} below on $K/F$ and $\Pi$ are satisfied. 

\section{Sorites on the sign}

\subsection{The notion of a good representation}

For a representation $\rho : G_K \rightarrow \Gl_n(L)$ that is a direct sum of absolutely irreducible and pairwise non isomorphic representations, and that satisfies (\ref{autobot}),
say that $\rho$ is {\it good} if for every irreducible factor of $\rho$ that appears with multiplicity one and satisfies (\ref{autobot}),  the sign of this factor is $+1$.

In this language, the theorem amounts to prove that $\rho_{\Pi,\lambda}$ is good,
which is good.
\subsection{Some trivial lemmas}
 
 In this paragraph, $\rho :  G_K \rightarrow \Gl_n(L)$ is a direct sum of absolutely irreducible representations, and satisfies (\ref{autobot}).

\begin{lemma} \label{twist}
If $\rho : G_K \rightarrow \Gl_n(L)$ is good, so is any twist by a character.
\end{lemma}
\begin{pf}
Note first that since $\rho$ satisfies (\ref{autobot}) for a character $\chi$, then if $\psi : G_K \rightarrow L^\ast$ is a character, 
$\rho\psi$ still satisfies (\ref{autobot}) for the character 
$\chi' = \chi \psi^{-1} \psi^\bot$ (which still satisfies $\chi'(g^c)=\chi'(g)$) 
and is a sum of absolutely irreducible pairwise non isomorphic factors,
namely the $\rho_i \psi$ where the $\rho_i$ are the factors of $\rho$.
Now if $\rho_i$ is an irreducible factor that satisfies (\ref{autobot}), 
a matrix $A$  that satisfies (\ref{A}) for $\rho_i$ and $\chi$
satisfies also (\ref{A}) for  
$\rho_i \psi$ and $\chi'$, hence the sign of $\rho_i$ and $\rho_i \psi$ are the same.
\end{pf}

\begin{lemma}\label{sub} If $\rho: G_K \rightarrow \Gl_n(L)$ is good, and $\rho'$ is a sub-representation of $\rho$ that satisfies (\ref{autobot}),
then $\rho'$ is good too.
\end{lemma}
\begin{pf} That one is {\it really} trivial.
\end{pf}

\begin{lemma}\label{rest} Let $F'$ be a totally real extension of $F$, and $K'=KF$. If $\rho|{G_{K'}}$ has the same number of irreducible components as $\rho$, and if those components are absolutely
irreducible, then $\rho|{G_{K'}}$ is good if and only if $\rho$ is good.
\end{lemma}
\begin{pf} If $\rho_i$ is an (absolutely) irreducible factor of $\rho$ that satisfies 
(\ref{autobot}), then $\rho_i{|G_{K'}}$ is still absolutely irreducible by hypothesis, still satisfies (\ref{autobot}),  and has obviously the same sign as $\rho_i$. The lemma follows.
\end{pf}

 \subsection{A specialization result} In this paragraph, $\anneau$ is a 
henselian discrete valuation domain with fraction field $L$ and residue 
field $k$, such that $2 \in \anneau^*$. We set also $G=G_K$ and assume 
that the character $\chi : G_K \rightarrow L^*$ actually falls into 
$\anneau^*$, thus it makes sense to talk about condition (\ref{autobot}) 
for $k$ or $L$-valued representations of $G$ (by a slight abuse of language,
we shall also denote by $\chi$ the residual character $G_K \rightarrow
k^*$). A simple but crucial 
observation for of our proof is the following:

\begin{prop} \label{spec} Assume that $\rho : G \rightarrow \Gl_n(\anneau)$ is such that 
$\rho \otimes L$ and $\rhob^\ses$ are a sum of absolutely irreducible pairwise non isomorphic representations
 and satisfy (\ref{autobot}). Then if $\rho \otimes L$ is good, so is
$\rhob^\ses$. 

Moreover, the converse holds if $\rhob^\ses$ has the same number of
irreducible factors as $\rho \otimes L$.
 \end{prop} 

Of course, in this statement $\rhob^\ses$ denotes the semisimplification of
the reduction $\rhob:=\rho \otimes_\anneau k$ of $\rho$. 

\begin{pf} Let $\rhob_1$ be a factor of $\rhob^\ses$ that
satisfies (\ref{autobot}).
Let $\tau_1,\dots,\tau_k$ be the irreducible factors of $\rho \otimes L$. For each of them
we can choose a stable $\anneau$-lattice, and see them as representations of $G$
over $\anneau$.
We have $\rhob^\ses = \oplus_{i=1}^k \overline{\tau_i}^\ses$ so $\rhob_1$
appears in exactly one
of the $\overline{\tau_i}^\ses$, say    $\overline{\tau_1}^\ses$. Moreover, 
$\overline{\tau_1}^\bot$ is isomorphic to $\overline{\tau_i} \chi$ for some
$i \in \{1,\dots,k\}$.
But it follows that $\overline{\tau_i}^\ses$ contains $\rhob_1$ (since
$\rhob_1$ satisfies
(\ref{autobot})), so the only possibility is that $i=1$. In other words,  
 $\tau_1$ satisfies (\ref{autobot}), and replacing $\rho$ by $\tau_1$,   
 we are reduced to prove the lemma with the supplementary assumption that
 $\rho \otimes L$ is absolutely irreducible. In that case, the proposition is \cite[Lemma 1.8.8]{BC}.
 
  Since this proposition is really one of the main tool used in our proof, 
and since the proof of \cite[Lemma 1.8.8]{BC} is a little bit difficult to separate 
from the  other concerns of \cite[\S1.8]{BC}, let us sketch it here for 
the convenience of the reader, trying to be as pedagogical as possible. 

Note first that the basic point that makes 
the result not obvious is that there is no reason so that we can find a matrix 
$A$ for $\rho \otimes L$ as in (\ref{A}) with $A \in \Gl_n(\anneau)$. A priori we just have $A \in \Gl_n(L)$, and it is therefore not possible to reduce 
$(\ref{A})$ mod $\m$.   

There is one case, however, where a simple proof is possible: assume that
$\rhob^\ses$ is absolutely irreducible. In this case, the representations
 $\rho^\bot$ and $\rho \chi$ over $\anneau$, being isomorphic over $L$ and 
residually absolutely irreducible, are isomorphic over $\anneau$ by a theorem of Serre and Carayol. In other words, we can find 
a matrix $A \in \Gl_n(\anneau)$ such that $(\ref{A})$ holds, and reducing this modulo $\m$, we get that $\rho \otimes L$ and $\rhob$ have the same sign in this case. Note that this proves also the last assertion of the theorem (in all cases!).

The proof of the general case consists in reducing to the residually irreducible case. This is 
not possible, however, if we keep working with representations of groups only. 
We have to work in the larger world of representations of algebras instead. As we saw, we may assume that $\rho \otimes L$ is absolutely irreducible, and we set $\rhob^\ses = \oplus_i \rhob_i$.

Let $R$ be the algebra $\anneau[G]$, and $S = \rho(R) \subset M_n(\anneau)$. We have $S\otimes_\anneau L=M_n(L)$. The algebra $S$
is provided with a natural $\anneau$-algebra anti-automorphism $\tau$, induced by the one on $R$
defined on $g \in G$ by $g \mapsto \chi(g)^{-1} (g^c)^{-1}$. Explicitly, by (\ref{A}), we have for $M \in S$, 
\begin{equation}\label{formuletau} \tau(M) = {}^t A^{-1} {}^t M {}^t A,\end{equation} and by our sign assumption ${}^t A = A$ : the involution $\tau$ is a symmetric involution of the matrix algebra $S\otimes_\anneau L = M_n(L)$.
\newcommand{\Sb}{\overline{S}}

On the other hand, let $\Sb$ denote the image of $k[G]$ in the $k$-endomorphisms of the representation $\rhob^\ses=\oplus_i \rhob_i$. Then $\Sb \simeq \prod_i M_{n_i}(k)$ ($n_i=\dim \rhob_i$, $\sum_i n_i=n$) and 
$\Sb$ is also provided with a natural $k$-algebra anti-automorphism $\tau$ as above. Moreover, there is a natural surjective $\anneau$-algebra map $S \rightarrow \Sb$
which is $\tau$-equivariant.

Let us denote by $\epsilon_i \in \Sb$ be the central idempotent corresponding to $\rhob_i$. It is well known that $\epsilon_i$ can be lifted as an idempotent $e_i$ of $S$ as $\anneau$ is henselian and $S$ finite over $\anneau$. However, we need a  more precise lifting result. Let us fix an $i$ such that $\rhob_i$ satisfies (\ref{autobot}), then we have $\tau(\epsilon_i)=\epsilon_i$. What we need is an idempotent $e_i$ in $S$ lifting $\epsilon_i$, such that $\tau(e_i)=e_i$. The existence of such 
an idempotent is easy to prove: first choose any lift $x \in S$ of 
$\epsilon_i$ and let $S_0$ be the sub-$\anneau$-algebra generated by 
$\frac{1}{2}(x+\tau(x))$. Obviously, $\tau$ fixes any element of $S_0$.
The restriction of the natural surjection 
$S \rightarrow \Sb$ to $S_0$ is onto a $k$-subalgebra
of $\Sb$ that contains the image of $\frac{1}{2}(x+\tau(x))$, that is 
$\epsilon_i$. Thus, defining $e_i$ as a lift of $\epsilon_i$ in $S_0$ does 
the job. (This result is the trivial case of \cite[Lemma 1.8.2]{BC}.) As $\rhob_i$ is absolutely irreducible and has multiplicity one in $\rhob^\ses$ it actually turns out that the rank of $e_i$ is 
$n_i=\dim \rhob_i$, and that $e_iSe_i \simeq M_{n_i}(\anneau)$. Replacing $\rho$ by a conjugate if necessary, we may then assume that $e_i$ is a diagonal idempotent of rank $n_i$ in $M_n(L)$.

Applying (\ref{formuletau}) to $M=e_i$ we get $Ae_i = {}^t e_i {}^t A$, that is $A e_i$ is symmetric. In other words, $\tau$ induces a symmetric involution on $e_iSe_i \simeq M_{n_i}(\anneau)$.
As a consequence, $\tau$ also induces a symmetric involution on $\epsilon_i\Sb\varepsilon_i = \End_k(\rhob_i)$, which exactly means that the sign of $\rhob_i$ is $+1$, QED. \end{pf}


\section{Proof of the main theorem} 

\subsection{Proof of theorem~\ref{main} under special hypotheses} \label{proofone}

We shall first prove the theorem under a set of additional hypotheses
on the CM extension $K/F$, the automorphic representation $\Pi$ and the place $\lambda$. 

Let us call $p$ the residual characteristic of $\lambda$. Recall that the automorphic representation $\Pi$ defines an embedding $\iota : E(\Pi) \rightarrow \C$.
We fix once and for all algebraic closures $\Qb$ and $\Qpb$ of $\Q$ and $\Q_p$, as well as some embeddings $\iota_\infty :
{\overline{\Q}} \rightarrow \C$ and $\iota_p:
\overline{\Q}_p \rightarrow \C$ such that the induced map $\iota_p\iota_{\infty}^{-1}\iota
: E(\Pi)
\rightarrow {\overline{\Q}}_p$ factors through $E(\Pi)_\lambda$.

\medskip
\subsubsection{Some special hypotheses}

\begin{itemize}
\item[(H1)] Special Hypotheses 2.2 of \cite{harris}, that is
\begin{itemize}

\item[(H1a)] $K/F$ is unramified at all finite places
\item[(H1b)] $\Pi_v$ is spherical at all non-split non-archimedean places $v$ of $K$
\item[(H1c)] The degree $[F:\Q]$ is even.
\item[(H1d)] All primes of small residue characteristic relative to $n$ are split in $K/F$.
\end{itemize} 
\item[(H2)] Special Hypothesis 2.3 of \cite{harris}, that is 
for at least one real place $\sigma$ of $K$, the infinitesimal character of $\Pi_\sigma$ is sufficiently far from the walls.\footnote{Here and everywhere in this paper, and as in \cite{harris}, this will mean that the extremal weight of the associated algebraic representation of $\GL_n(\C)$ does not belong to a wall. For our purposes, we could even assume here that this holds for all archimedean places.} 
\item[(H3)] There\footnote{See \S\ref{generalnotations} for the notations used in this assumption.} is a place $v$ above $p$ in $F$ that splits in $K$, and
for $w$ a place of $K$ above $v$, $\Pi_w$ is unramified. Denote by 
$\{\varphi_1,\dots,\varphi_n\}$ the eigenvalues of
$L(\Pi_w|.|^{(1-n)/2})$. Then the Hodge-Tate weights $\{k_{i,\sigma}\}$ of $\rho_{\Pi,\lambda}|G_w$ and the slopes $\v(\varphi_j)$ are in sufficiently general position in the following sense: if $$c = \max_{i \in \{1,\dots,n\}}
\min_{j \in \{1,\dots,n\}} |\v(\varphi_i)-k_j|$$ then for all subsets $I$ and $J$
of $\{1,...,n\} \times \Hom(K_w,\Qpb)$ with $|I|=|J|<nd$, we have
$$|k_I - k_J| > (n+1)\cdot c.$$\end{itemize}
In (H3) above, $\v : \Qpb \rightarrow \overline{\Q}$ is the valuation such that $\v(p)$ is the ramification index of $p$ in $K_w$.

\subsubsection{The theorem}
We want to prove :
\begin{theorem} \label{H123} With the supplementary hypotheses (H1), (H2) and (H3),
Theorem~\ref{main} holds.
\end{theorem}
The rest of this subsection is entirely devoted to the proof of this theorem.

\subsubsection{Descent and transfer}

Let $m=n$ if $n$ is odd, and $m=n+1$ if $n$ is even, so that $m$ is always odd.
Let us call $\U(m)$ be the unitary group over $F$ attached to $K$ in $m$ variables
that is quasi-split at every finite places of $F$ and compact at every infinite place.
Since $m$ is odd, such a group always exists (uniquely up to isomorphism).
Actually, $\U(m)$ is simply the standard unitary group attached to the hermitian form
$\sum_{i=1}^m N_{K/F}(z_i)$ on $K^m$ (\cite[\S6.2.2]{BC}).

If $n$ is odd, that is if $n=m$, by hypothesis (H1) and Labesse's base change theorem \cite[chap. 4 book 1]{book}, we can
descend $\Pi$ to a representation $\pi$ of $\U(m)$ with $\pi_v \simeq \Pi_w$
for every place $w$ of $K$ split over $v$ in $F$ (with the natural identification
$\U(n)(F_v) \simeq \GL_n(K_w)$), and such that for each complex place $w$ of $K$ above a real place
$v$ of $F$, $\pi_v$ has the same infinitesimal character as
$\Pi_w$ (under the natural identification $\U(n)(K_w)\simeq
\GL_n(\C)$).

If $n$ is even, we use a result of endoscopic transfer due to Clozel, Harris, and Labesse
 (see \cite{book} and in particular \cite{harris}). 
Note first that using $\iota_\infty\iota_p^{-1}$, if $v=ww^c$ is as (H3) we may identify $\Hom(K_w,\Qpb)=\Hom(F_v,\Qpb)$ with subsets $\Sigma_v$ and $\Sigma_w$ of $\Hom(F,\R)$ and $\Hom(K,\C)$. Let us first fix
$$\mu : K^* \backslash \A_K^* \rightarrow \C^*$$

\noindent a Hecke character such that $\mu(c(x))^{-1}=\mu(x)$, and such that for each
$s \in \Sigma_v$, $\mu_s(z)=(\sigma_s(z)/\overline{\sigma_s(z)})^{\frac{1}{2}}$ where $\sigma_s \in \Sigma_w$ is associated to $s$ as above. This last assumption implies
that $\mu|\A_F^*$ coincides with the sign of $K/F$, and that $\mu$ does not
come by base change from a Hecke-character of $\U(1)$ (see e.g. \cite[\S
6.9.2]{BC}). Such a Hecke character always exists, and as $K/F$ is
unramified at all finite places, we can even assume (and we will) that it is
unramified at the finite places of $K$ which are either above $p$ or not
split above $F$. Let us choose another Hecke character $$\chi : K^* \backslash
\A_K^* \rightarrow \C^*$$
such that $\chi(c(z))^{-1}=\chi(z)$ but assume now that $\chi$ descends to
$U(1)$, i.e. that for each real place $s \in \Sigma_v$, $\chi_s(z)=\sigma_s(z/c(z))^{-a_\sigma}$ for some $a_s \in \Z$. Assume
also that $\chi$ is unramified at the finite places of $K$ which do not split
over $F$. Under hypotheses (H1) and (H2), and if all the $|a_s|$ are
big enough, by the aforementioned results of Clozel, Harris and Labesse, we
can transfer $\Pi$ to an automorphic representation $\pi$ of $\U(m)$ in such a way
that at every place $w$ of $K$ split over a place $v$ in $F$, we have
\begin{eqnarray} \label{lpv} L(\pi_v)=L(\Pi_w \mu_w) \oplus L(\chi_v).\end{eqnarray}
 Moreover, for each real place $v$ of $F$ and each complex place $w$ of $K$
above $w$, the infinitesimal character of $\pi_v$ is obtained from the one
of $\Pi_w \mu_w$ in the obvious way : in terms of the associated
Harish-Chandra's cocharacter, it is the direct sum of the one of $\Pi_w \mu_w$
and the one of $\chi_w$. 

In both cases ($n$ even or odd), the aforementioned authors actually
construct a $\pi$ which is moreover unramified at all the finite places of $K$ which are
not split over $F$ (we don't really need this, but this fixes ideas). 

\subsubsection{Consequences of (H3)}

When $n=m$ is odd, we set $\rho_\pi:= \rho_{\Pi,\Lambda}$
When $n$ is even, the $G_K$ representation of dimension $m$ attached to $\pi$ is by 
definition $$\rho_\pi := \rho_{\Pi,\lambda}(\mu|.|^{-\frac{1}{2}})\oplus \chi|.|^{(1-m)/2}.$$
Note that $\mu|.|^{-\frac{1}{2}}$ and $\chi|.|^{(1-m)/2}$ are both algebraic Hecke characters of $K$, and we identify them here with their $p$-adic
realization given by $\iota_\infty$ and $\iota_p$. By assumption, $\mu|.|^{-\frac{1}{2}}$ is actually unramified at the place $w$, and $\chi|.|^{(1-m)/2}$ is crystalline, and we shall denote by $\varphi_\mu$ and $\varphi_\chi \in \Qpb^*$ their associated Frobenius eigenvalue. 

If $n$ is even, so $m=n+1$, we set for each $\sigma \in \Hom(K_w,\Qpb)$, 
$$k_{m,\sigma}:=\frac{m-1}{2}+a_{\sigma}$$ (where $\sigma$ is viewed as an element of $\Hom(F,\R)$ as above). Thus, in any case, the $k_{i,\sigma}$ for $i=1,\dots,m$ and $\sigma \in \Hom(K_w,\Qpb)$ are the Hodge-Tate weights of $\rho_\pi|G_w$.
We shall use for this extended collection $\{k_{i,\sigma}\}$ with all $i \in\{1,\dots,m\}$, and for a subset $I$ of $\{1,\dots,m\} \times \Hom(K_w,\Qpb)$, the notation $k_I$ analogous to the one in \S\ref{generalnotations}. 

If $n$ is even, we set $\varphi'_i:=\varphi_i\varphi_\mu$ for $i<m$ and 
$\varphi'_m:=\varphi_\chi$. We have $\v(\varphi'_i)=\v(\varphi_i)$ for $i<m$ and 
$\v(\varphi'_m)=k_m$. When $n=m$ is odd, we shall simply set $\varphi'_i:=\varphi_i$.
Thus, in both cases, $\varphi'_1,\dots,\varphi'_m$ are the Frobenius eigenvalues of $L(\pi_w|.|^{(1-m)/2})$.

If $n$ is even, we precise now our choice of $\chi$. We assume that 
$k_{\sigma, m} = a_{\sigma}+\frac{d(m-1)}{2}$ are all big with respect to
the $k_i$ and $\v(\varphi'_i)$ for $i=1,\dots,n$, and also that they 
are set sufficiently far  apart so that any non trivial sum of the form 
$\sum_{\sigma \in S} \pm k_{\sigma,m}$, where $S \subset \Hom(K_w,\Qpb)$, 
is big in the same sense as above. This is of course always possible.
With those assumptions:

\begin{lemma}\label{H3prime} \begin{itemize}\item[(i)] The representation $\pi_v$ is a fully induced unramified principal series, and the eigenvalues of $L(\pi_w|.|^{(1-m)/2})$ are $\varphi'_1,\dots,\varphi'_m$. 
\item[(ii)] We  have $c = \max_{i\in\{1,\dots,m\}}
\min_{j\in\{1,\dots,m\}} |\v(\varphi'_i)-k_j|$, and for all distinct subsets $I$ and $J$ of $\{1,\dots,m\} \times {\rm Hom}(K_w,\Qpb)$ and with $|I|=|J|<md$, we have
$|k_I - k_J| > m\cdot c$.
\end{itemize}
\end{lemma}
\begin{pf}
By (H3) and, if $n$ is even, by (\ref{lpv}), $\pi_v$ is unramified. Moreover, the eigenvalues of $L(\pi_w|.|^{(1-m)/2})$ are $\varphi'_1,\dots,\varphi'_m$, 
and no quotient of those eigenvalues is equal to $q$, the cardinal of the residue field of
$K_w$. Indeed, a well-known
result of Jacquet-Shalika asserts that for $i=1 \dots n$, the complex numbers
$q^{(1-m)/2}\iota_\infty\iota_p^{-1}(\varphi'_i)$ are $< q^{1/2}$ in
absolute value, and $q^{(1-m)/2}\iota_\infty\iota_p^{-1}(\varphi'_m)$ has
norm $1$ by construction. Hence $\pi_v$ is a full unramified principal series
by Zelevinski's theorem, which is (i).

For (ii), there is nothing to prove if $n=m$ is odd. Assume $n$ even so $m=n+1$. 
Let us note that for $i=m$, we have 
$\min_{j\in\{1,\dots,m\}} |\v(\varphi'_i)-k_j| = 0$ since $\v(\varphi'_m)=k_m$.
For $i \leq n$, the minimum $\min_{j\in\{1,\dots,m\}} |\v(\varphi'_i)-k_j| = 0$ is
not realized for $j=m$ because  $k_m$ is much too big.

Hence \begin{eqnarray*} \max_{i\in\{1,\dots,m\}}
\min_{j\in\{1,\dots,m\}} |\v(\varphi'_i)-k_j| &=& \max_{i\in\{1,\dots,n\}}
\min_{j\in\{1,\dots,n\}} |\v(\varphi'_i)-k_j|\\ &=& c.\end{eqnarray*}

It remains to prove that $|k_I-k_J| > mc = (n+1)c$. Let $I_0$ (resp. $J_0$)
be the subset of $I$ (resp. of $J$) of pairs $(i,\sigma)$ with $i=m$.
 If $I_0=J_0$, then $k_I - k_J = k_{I-I_0} - k_{J-J_0}$ and since $I-I_0$, $J-J_0$, are distinct subsets of same cardinality of $\{1,\dots,n\} \times \Hom(K_w,\Qpb)$, the desired inequality comes directly from (H3). If $I_0 \neq J_0$, $k_I-k_J$ contains, in
addition of a bounded numbers of terms $\pm k_{i,\sigma}$ for $i \leq n$, a non trivial sum of the form 
$\sum_S \pm k_{\sigma,m}$, where $S \subset \Hom(K_w,\Qpb)$, hence $|k_I-k_J|$ is again greater than $mc$.
\end{pf}

\subsubsection{Eigenvarieties and their families of Galois representations}\label{pareigen}

We are ready now to start the deformation argument. Let 
$U= \prod_v U_v \subset \U(n)(\A_{F,f})$ be a compact open subgroup 
such that $\pi^U \neq 0$, and assume that $U_v$ is a Iwahori subgroup 
for the place $v$ of (H3), and that $U_v$ is hyperspecial for all places of $K$ 
that are not split over $F$.

From now on, we shall reserve the notation $v$ for the place of $F$ of hypothesis (H3), 
and $w$ for one of the place of $K$ above $v$. We shall denote by $d$ the degree of the field $F_v = K_w$
over $\Q_p$. To $U$, the place $v$ and $(\iota_p,\iota_\infty)$, we can attach an
eigenvariety $X=X_{U,v,(\iota_\infty,\iota_p)}$ for the group $\U(m)/F$ which is a reduced rigid analytic space over $\Q_p$ of equidimension\footnote{It is not necessary here to let the weights corresponding to the other possible places of $F$ above $p$ move, but we could have, and the eigenvariety would then have dimension $n[K:\Q]$.} $md$. The construction of this eigenvariety follows essentially verbatim from the method of \cite{Ch} 
(which is only written in
the case $F=\Q$, but see e.g. \cite{buzzard} for the setting for a general $F$ 
in dimension\footnote{The situation here is actually even simpler because the center $\U(1)$ of $\U(m)$ is anistropic over $\R$.} $2$. Details should appear elsewhere as part of the book project~\cite{book}). Alternatively, we
may deduce it from the work of Emerton in \cite[\S3.2]{emerton}. 

		By Labesse's base change theorem, if $\pi'$ is any automorphic representation of $\U(m)$ which is unramified outside the split finite places of $K/F$, then $\pi'$ admits a strong base change to $\GL_m/K$, hence a Galois representation by Expected Theorem~\ref{galois}. As explained in \cite[Chap. 7.5]{BC} (or \cite{Ch}), this is enough to equip $X$ with a continuous $m$-dimensional pseudocharacter $T : G_K \rightarrow \OO(X)$ of dimension $m$. The eigenvariety $X$ and this $T$ satisfy a number of properties and we will only list below the ones we shall need. If $x \in X(\bar \Q_p)$ we note $T_x$ the evaluation of $T$ at $x$ and $\rho_x$ the semi-simple representation $G_K \rightarrow \Gl_m(\bar \Q_p)$ of trace $T_x$. There is :
		
\begin{itemize}
\item[(i)] Zariski-dense and accumulation subsets $Z^{\rm reg} \subset Z \subset X(\Qpb)$ of {\it classical points},
\item[(ii)] a set of $dm$ analytic functions\footnote{Again, we shall use for this collection $\{\kappa_{i,\sigma}\}$, and for a subset $I$ of 
$\{1,\dots,m\}\times \Hom(K_w,\Qpb)$ (resp. an $i \in \{1,\dots,m\}$), the notation $\kappa_I$ (resp. $\kappa_i$) analogous to the one in \S\ref{generalnotations}.} $\kappa_{1,\sigma},\dots,\kappa_{m,\sigma}$ where
$\sigma$ runs over the embeddings $K_w \rightarrow \Qpb$,
\item[(iii)] a set of locally constant functions $s_1,\dots,s_m : X(\bar \Q_p) \rightarrow \Q$,
\end{itemize} 
satisfying the following conditions :
\begin{itemize}
\item[(a)] if $z \in Z$, $\rho_z|G_w$ is crystalline.
\item[(b)] if $z \in Z$, the ordered Hodge-Tate weights of $\rho_z|G_w$ are $\{\kappa_{i,\sigma}\}$
\item[(c)] let $C$ be any real number and $Z_C := \{z \in Z^{\rm reg}, |\kappa_I(z)-\kappa_J(z)| >
C$ for all 
distinct subsets
$I,J$ of $\{1,\dots,m\} \times \Hom(K_w,\Qpb)$ and $|I|=|J|<md\}$. Then $Z_C$ is Zariski dense and accumulation in $X$. 
\end{itemize}

Moreover, the classical points $z$ in $Z$ correspond to pairs $(\pi(z),\RR(z))$ where 
$\pi(z)$ is an automorphic representation of $\U(m)$ such that $\pi(z)^U \neq 0$ and $\RR(z)=(\varphi_1,\dots,\varphi_m)$ is an accessible refinement\footnote{Recall that an {\it refinement} of an irreducible smooth representation $\rho$ of $\GL_m(K_w)$ such that $\rho^I \neq 0$ for $I$ an Iwahori subgroup is an ordering $(\varphi_1,\dots,\varphi_m)$ of the eigenvalues of $L(\rho)(\Frob_w)$. It is said {\it accessible} if $\rho$ appears as a {\it sub-representation} of the induced representation ${\rm Ind}_B^{\GL_m(K_w)} \chi \delta_B^{1/2}$ where $B$ is (say) the upper Borel subgroup, $\delta_B$ is modulus character, and $\chi$ the (unramified) character of the diagonal torus sending $(x_1,\cdots,x_m)$ to $\prod_{i=1}^m \varphi_i^{\v(x_i)}$ (see \cite[\S6.4.4]{BC}).}
of $\pi(z)_w|.|^{(1-m)/2}$, in the following sense :
$\rho_z$ is the Galois representation attached to the base change of $\pi(z)$ to $\Gl_m/K$ by Expected
Theorem~\ref{galois} and for each $i=1,\dots,m$, $\v(\varphi_i)=s_i(z)+\kappa_i(z)$. 

\begin{itemize}\item[(d)] If $z \in Z$ parameterizes $(\pi(z),\RR(z)=(\varphi_1,\dots,\varphi_m))$, then for all $i$ we have $\v(\varphi_i)=s_i(z)+\kappa_i(z)$.
\end{itemize}

The subset $Z^{\rm reg} \subset Z$ parameterizes refined automorphic representations $(\pi,\RR)$ satisfying some additional properties, and for our concerns here we shall simply assume that they are those $(\pi,\RR)$ such that $\pi_v$ unramified and such that for each real place $s$, the infinitesimal character of $\pi_s$ is sufficiently far from the walls. Under this latter condition, the base change of an automorphic representation of $\U(m)$ is not necessarily cuspidal, but always associated to a decomposition $m_1+\dots+m_r=m$ and a $r$-tuple $(\pi_1,\dots,\pi_r)$ of cuspidal (polarized, cohomological) automorphic representations $\pi_i$ of $\GL_{m_i}(\AAA_K)$ ; moreover each $\pi_i$ satisfies property (H2) in dimension $m_i$ and is unramified at $v$. In particular, for a $z \in Z^{\rm reg}$, the characteristic polynomial of the crystalline Frobenius of $\rho_z|G_w$ coincides with the polynomial $P_w(T)$ associated to $\pi_w|.|^{(1-m)/2}$ by Expected Theorem~\ref{galois} (see \S\ref{generalnotations}), and we also have the following :

\begin{itemize}\item[(d')] If $z \in Z^{\rm reg}$, then the $m$ slopes of the crystalline Frobenius of $\rho_z|G_w$ are the $s_i(z)+\kappa_i(z)$, fo $i=1,\dots,m$.
\end{itemize}

\subsubsection{Choice of a refinement}

Going back to the representation $\pi$ introduced above, if we choose an accessible refinement $\RR$ of $\pi_v$,
there is a point $z_0 \in Z$ corresponding to $(\pi,\RR)$.

\begin{lemma} There exists a refinement $\RR$ of $\pi_v$ such that the pseudocharacter $T$ is generically irreducible in a neighborhood of the corresponding
point $z_0$.
\end{lemma}
\begin{pf} We shall eventually show that the conclusion holds for $T|G_w$. Note that by construction, for all $\sigma \in \Hom(K_w,\Qpb)$ and $i \in\{1,\dots,m\}$,  
$\kappa_{i,\sigma}(z_0)=k_{i,\sigma}$. Let us first renumber the 
$\varphi'_i \in \Qpb^*$ so 
that $|\v(\varphi'_i) - k_i| = \min_j |\v(\varphi'_i)-k_j|$. By Lemma~\ref{H3prime} (ii) there is one and only one way to do so, and this being done we have 
$\v(\varphi'_1) < \v(\varphi'_2) < \dots < \v(\varphi'_m)$ (strict inequalities). Then consider a transitive permutation $\sigma$ of $\{1,\dots,m\}$.
We choose the refinement $$\RR=(\varphi'_{\sigma(1)},\dots,\varphi'_{\sigma(m)}).$$
Since $\pi_v$ is a full unramified principal series by Lemma~\ref{H3prime} (i), all the refinements of $\pi_v$ are
accessible, so $\pi$ together with $\RR$ defines a point $z_0$.

Before proving the irreducibility property of the lemma, let us observe a combinatorial property of this refinement. We have by 
definition $\kappa_i(z_0)=k_i$ and
$s_i(z_0) = \v(\varphi'_{\sigma(i)}) - k_i$. We claim that for any non-empty proper subset 
$I \subset \{1,\dots,m\}$, \begin{eqnarray}\label{comb} \sum_{i \in I} s_i(z_0) \neq 0. \end{eqnarray}
  Indeed, we compute \begin{eqnarray*}
  |\sum_{i \in I} s_i(z_0)| &=& |\sum_{j \in J} \v(\varphi'_j) - \sum_{i \in I} k_i| \text{ where $J=\sigma(I)$} \\
  &=& |(\sum_{j\in J} k_j - \sum_{i \in I} k_i) + \sum_{j \in J} (\v(\varphi'_j)-k_j) | \\
  &\geq&  |(\sum_{j\in J} k_j - \sum_{i \in I} k_i)| - \sum_{j \in J} | \v(\varphi'_j)-k_j | \\
  &>& mc - mc \text{ by Lemma~\ref{H3prime}(ii) as $I\neq J$} \\
  &=& 0.
  \end{eqnarray*}

Let us choose now some affinoid neighborhood $\Omega$ of $z_0 \in X$ on which the $s_i$ are constant and in which $Z_C$ is Zariski-dense for $C = \sum_{i=1}^m  |s_i(z_0)|$. We claim that for
every point $z$ of $Z_C \cap \Omega$, $\rho_z|G_w$ is irreducible. Indeed if it was not, it would have a sub-representation of dimension $0<r<m$, and by the weak admissibility of $D_\crys(\rho_z|G_w)$ there would exist a subset $I \subset \{1,\dots,m\}$ of cardinal $r$, and a subset $J \subset \{1,\dots,m\} \times \Hom(K_w,\Qpb)$ with $|J|=rd$, such that
  $$ \sum_{i \in I} (\kappa_i(z)+  s_i(z)) = \kappa_J(z).$$
(here we use that $z \in Z^{\rm reg}$ and property (d') of eigenvarieties.) Since $z \in Z_C$, we see at once that $I \times \Hom(K_w,\Qpb)=J$. But this implies 
  $$0 = \sum_{i \in I} s_i (z),$$
  a contradiction with (\ref{comb}) as $s_i(z)=s_i(z_0)$ for all $i$. 
\end{pf}

 \subsubsection{End of the proof}

Let $\Omega \subset X$ be the neighborghood defined above of the point $z_0$, and let $A$ be a complete discrete valuation ring, with a map of $\spec A$ to the spectrum of the rigid local ring
$\anneau_{z_0}$ of $X$ at $z_0$ which sends the special point of $\spec A$ to $z_0$ and the generic
point to the generic point of any irreducible component of $\Omega$ containing $z_0$. Let us call $L$ the fraction field of $A$ and $\m$ its maximal ideal. By pulling back the pseudocharacter $T$ over $A$, we get a representation
 $\rho : G_K \rightarrow \GL_m(A)$ such that $\rho \otimes L$ is absolutely irreducible and satisfies (\ref{autobot}) (for $\chi=\Q_p(m-1)$) and $$\rhob^\ses = \left\{ \begin{array}{ll}  \rho_{\Pi,\lambda} \,\, \, \text{ if \,\,\,}m=n,\\ \rho_{\Pi,\lambda}(\mu|.|^{-1/2}) \oplus (\chi|.|^{(1-m)/2})\,\,\, \text{if \,\,\,}\, m=n+1.\end{array}\right.$$ 
Since $\rho \otimes L$ is absolutely irreducible and satisfies  (\ref{autobot}) it has a sign that can only be $+1$. Hence it is good, and so is $\rhob^\ses$ by Prop.~\ref{spec}, hence $\rho_{\Pi,\lambda}$ by Lemmas~\ref{sub} and~\ref{twist}, QED.

\subsection{Weakening of the hypothesis (H3), removal of (H2)}

We consider the following hypothesis which is a weakening of (H3)
\begin{itemize}
\item[(H4)] 
There is a place $v$ above $p$ in $F$ that splits in $K$, and for $w$ a place of $K$ above $v$, $\Pi_w$ has a non-zero vector invariant by a Iwahori subgroup of $\Gl_n(K_w)$.
\end{itemize}

\begin{theorem} \label{H14}
With the supplementary hypotheses (H1), and (H4),
Theorem~\ref{main} holds.
 \end{theorem}
 
 We shall argue by induction on $n\geq 1$, there is nothing to show if $n=1$.
 
 Let $\U(n)$ be the $n$-variables unitary group over $F$ attached to $K/F$ that is
 quasi-split at every finite place and compact at every infinite place. Hasse's principle shows that this group exists, even if $n$ is even, by condition (H1c) (see e.g. \cite[Lemma 3.1]{harris}). Moreover, condition (H1) and Labesse's base change theorem also ensures that $\Pi$ descends to an automorphic representation $\pi$ for $\U(n)$. Again, $\pi$ is unramified at non split finite places of $K/F$  and for each complex place $w$ of $K$ above a real place
$v$ of $F$, $\pi_v$ has the same infinitesimal character as
$\Pi_w$ (under the natural identification $\U(n)(K_w)\simeq
\GL_n(\C)$).

 Let $U= \prod_v U_v \subset \U(n)(\A_{F,f})$ be a compact open subgroup such that $\pi^U \neq 0$, and assume that $U_v=I_v$ for the place $v$ of (H4), and that $U_v$ is hyperspecial for each place $v$ of $K$ that is not split over $F$. Let $X$ be the  eigenvariety attached to $U$, $v$ and $\iota_\infty, \iota_p$. Of course, all what we said for the eigenvarieties of $\U(m)$ with $m$ odd also applies verbatim to this $X$ by replacing the letter $m$ by letter $n$, and we shall not repeat it.

 	Let $z_0$ be the point corresponding to $\pi$ together with some accessible refinement of $\pi_v$. As a general fact, there is always such a refinement and we choose it anyhow here.  
 
  Let $c$ be the maximum of the $|s_i(z_0)|$ and $C > nc$. Let $\Omega \subset X$ be an open affinoid of $X$ containing $z_0$, in which $Z_C$ is Zariski-dense, and over which the $s_i$ are constant. We claim that for $z \in Z_C$, $\rho_z$ is good. Indeed, let $\Pi(z)$ be Labesse's base change of $\pi(z)$ to $\Gl_n(\A_K)$. As $z \in Z^{\rm reg}$, and as explained in \S\ref{pareigen}, there exists a decomposition $n=n_1+\dots+n_r$ and cuspidal automorphic representations $\Pi_i$ of $\GL_{n_i}(\AAA_K)$, satisfying (H1b), (H2) and unramified at $v$, such that $$\rho_z = \bigoplus_{i=1}^r \rho_{\Pi_i,\lambda}\otimes \chi_i,$$
for some characters $\chi_i : G_K \rightarrow \Qpb^\ast$. If $r>1$, then $\rho_z$ is good by induction and Lemma~\ref{twist}. If $r=1$, then $\Pi(z)$ is cuspidal and it satisfies (H2) and (H3) by construction, so $\rho_z$ is good by Theorem~\ref{H123}. (To check (H3), remark that for such a $z$, and for each $i \in \{1,\dots,n\}$, we have ${\rm Min}_j (s_i(z)+k_i(z)-k_j(z))=s_i(z)=s_i(z_0)$.)   
 
		Let $W$ be any irreducible component of $\Omega$ containing $z_0$, and ${\rm Frac}(W)$ its associated function field. As $Z_C$ is Zariski-dense in $\Omega$ we may find a $z \in Z_C \cap W$ such that the pseudocharacters $T_z$ and $T \otimes_{\anneau(\Omega)}{\rm Frac}(W)$ have the same number of irreducible factors. Such factors are necessarily absolutely irreducible here, by \cite[Thm. 1.4.4 (iii)]{BC}. 
Arguing as in the preceding section, let $A$ be a complete discrete valuation ring with a map of $\spec A$ to 
$\spec \anneau_{z}$ which sends the special point of $\spec A$ to $z$ and its generic
point to ${\rm Frac}(W)$, and let $\rho : G_K \rightarrow \GL_n(A)$ be a representation with trace $T$ such that $\rho \otimes_A L$ is a direct sum of absolutely irreducible, distinct, representations (use e.g. \cite[Prop. 1.6.1]{BC}). As we saw,  $\rhob=\rho_z$ is good, hence so is $\rho \otimes_A L$ by Prop.\ref{spec}, as well as $\rho \otimes_A {\rm Frac}(W)$ for any irreducible component $W$ containing $z_0$. But arguing back now at the point $z_0$ as in the preceding section, we obtain that $\rho_{z_0}=\rho_{\Pi,\lambda}$ itself is good as a specialization of a good representation, and we are done.
  
\subsection{Removal of Hypotheses (H1) and (H4)}
 
We now prove Theorem~\ref{main}. 

\begin{lemma} Let $L$ be a finite extension of $\Q_p$ and $\rho : G_K \longrightarrow \GL_n(L)$ a continuous representation which is a direct sum of absolutely irreducible
representations. There is a finite Galois extension $M/K$ such that for every finite extension $K'/K$ linearly disjoint from $M$, $\rho$ and $\rho|{G_{K'}}$ have 
the same number of irreducible factors, and the irreducible factors 
of $\rho|G_{K'}$ are absolutely irreducible.
\end{lemma}
\begin{pf} We can assume without loss of generality that $\rho$ is absolutely 
 irreducible. In particular, there exists $n^2$ elements $g_1,\dots,g_{n^2}$ such that
 the $\rho(g_i)$'s generate $M_n(L)$ as a $L$-vector space. Since $G_K$ has a basis
 of neighbourhoods of $1$ that are open normal subgroups, and since $\rho$
and the determinant are continuous, there is an open normal 
subgroup $U$ of $G_K$ such that if for all $i=1,\dots,n^2$,
 $g'_i \in g_i U$, then the $\rho(g'_i)$'s still generate $M_n(L)$.
 Set $M=\bar K^U$, so $M$ is a finite Galois extension of $K$.
 
 If $K'$ is a finite extension of $K$ which is linearly disjoint from $M$, so is its Galois closure. Hence we may assume that $K'$ is Galois over $K$. Thus 
 $\Gal(K'M/K')$ is naturally isomorphic to $\Gal(K/M)$.
 For every $i$, choose $g'_i$ in $G_{K'}$ whose image in $\Gal(K'M/K')$ is sent to $g_i$
 by the above isomorphism. This implies that $g'_i g_i^{-1} \in U$, hence the 
 $\rho|G_{K'}(g'_i)$'s generate $M_n(L)$, and $\rho|G_{K'}$ is absolutely irreducible.
 \end{pf}
 
By \cite[Prop. 4.1.1 and Thm 4.2.2]{harris} for any finite extension $M/K$ 
there exists a totally real solvable Galois extension $F'/F$ such 
that $K'=KF'$ is linearly disjoint to $M$ and such that 
Arthur-Clozel's base change $\Pi_{K'}$ and $K'/F'$ satisfy hypotheses (H1) and (H4). We apply this to an $M$ as in the lemma above. By Theorem~\ref{H14} we know that $(\rho_{\Pi,\lambda})_{|G_K'}$ is good, and by the lemma above, that it has the same number of 
irreducible components as $\rho_{\Pi,\lambda}$.

Hence  by Lemma~\ref{rest} $\rho_{\Pi,\lambda}$ is good, QED.

\par \bigskip
\end{document}